\documentclass[12pt]{article}
\usepackage{a4,amsmath,amsthm,amssymb,amsfonts,hyperref,enumerate}

\theoremstyle{plain}
\newtheorem{thm}{Theorem}      
\newtheorem{cor}{Corollary}
\newtheorem{lem}{Lemma}

\theoremstyle{definition}
\newtheorem{rem}{Remark}  

\newtheorem{prob}{Problem} 
\newtheorem{quest}{Question}

\newcommand{\abs}[1]{{\left| {#1} \right|}}
\newcommand{\p}[1]{{\left( {#1} \right)}}
\newcommand{\Oh}[1]{{O \p{#1}}}

\author{Johan Andersson\footnote{Email:johana@math.su.se}}

\begin{document}

\title{On some power sum problems of Montgomery and Tur\'an}



\maketitle

\begin{abstract}
  We  use an estimate for character sums over finite fields of Katz to solve open problems of Montgomery and Tur\'an.  Let $h \geq 2$ be an integer. We prove that
    $\inf_{\abs{z_k}= 1} \max_{\nu=1,\ldots,n^h} \abs{\sum_{k=1}^n
      z_k^\nu} \leq (h-1+o(1)) \sqrt n.$
 This gives the right order of magnitude for the quantity and improves on a bound of Erd{\H os}-Renyi by  a factor of the order $\sqrt{\log n}$.
\end{abstract}

\section{Introduction}
One of the simplest and most striking problems in the Tur\'an power sum theory is to determine the quantity
\begin{gather*}
  (\star)= \inf_{\abs{z_k} = 1}\max_{\nu=1,\ldots,m} \abs{\sum_{k=1}^n    z_k^\nu} 
\end{gather*}
for various choices of integers $n,m$. The purpose of this paper is to find the correct order of magnitude for this quantity when $m \sim n^B$. For any fixed $B>1$ and $m=\lfloor n^B \rfloor$ we prove that $(\star) \asymp \sqrt n$. This solves open problems of Hugh Montgomery and Paul Tur\'an.

We first present some well known results about $(\star)$. By letting
\begin{gather*}
  z_k= e \p{\frac k {n}}, \qquad \qquad (e(x)=e^{2 \pi i x})
\end{gather*}
we see that $(\star)=0$ if $1 \leq m \leq n-1$. Tur\'an \cite{Turan2} proved that $(\star)=1$ if $m=n$. In general it is difficult to determine the quantity $(\star)$ exactly but in \cite{Andersson2} we proved that $(\star)=\sqrt{n-1}$  if $m=n^2-n$ and $n-1$ is a prime power, and if $m=n^2-j$ for $2 \leq j \leq n-1$ and $n$ is a prime power then $(\star)=\sqrt n$.

When exact values can not be determined we are interested in obtaining asymptotic estimates for $(\star)$ when $n,m \to \infty$.  In the special case $n^{1+\delta} < m \leq n^2$ we have proved in \cite{Andersson3}  that $(\star) \sim \sqrt n$.  

Finding asymptotic estimates can be quite difficult as well and the next question is to ask for choices of $m,n$ where we can  find the order of magnitude for $(\star)$. 
Montgomery proved (\cite{Montgomery}, page 100, Theorem 10) that 
 \begin{gather} \label{mont}
\sqrt{nB} \ll \max_{\nu=1,\ldots,n^B} \abs{\sum_{k=1}^n    z_k^\nu}
\end{gather}
uniformly for $1+\delta \leq B \leq n$ when $|z_k|=1$. Erd\H os-Renyi \cite{ErdosRenyi} have used probabilistic reasoning to show that there exists an $n$-tuple 
$(z_1,\ldots,z_n)$ of unimodular complex numbers such that
\begin{gather} \label{erd}
  \abs{\sum_{k=1}^n    z_k^\nu} \leq \sqrt{6 n \log(m+1)} \qquad (\nu=1,\ldots,m)
\end{gather}
thus giving a general upper bound for $(\star)$. Let us now consider the case $m=\lfloor n^B \rfloor$ for a fixed constant $B>1$.  Leenman-Tijdeman \cite{LeenmanTijdeman} have given an explicit construction which yields the same order of magnitude as Eq. \eqref{erd}.  Unfortunately these upper estimates differ by a essentially a $\sqrt{\log n}$ factor from the lower bound Eq. \eqref{mont}. 
This inspired Hugh Montgomery  to state the following problem  
\begin{prob} (Montgomery, \cite{Montgomery}, page 197, Problem 13)
 Show that for any positive $B$ there exist complex numbers $z_1,\ldots,z_n$ such that $|z_k|=1$ for all $k$ and
$$ \left|\sum_{k=1}^n z_k^\nu \right| \ll_B  \sqrt n. \qquad \qquad (\nu=1,\ldots,\lfloor n^B \rfloor)$$
\end{prob}
In contrast Tur\'an had previously  stated the following problem.
\begin{prob} (Tur\'an, \cite{Turan} page 197, Problem 54) Does there exist an $w(x) \nearrow \infty$ such that
\begin{gather*}
              b_j>0, \qquad \abs{z_j}=1, \qquad j=1,\ldots,n \\ \intertext{implies for $g(\nu)=\sum_{j=1}^n b_j z_j^\nu$ the inequality}
      \max_{1 \leq \nu \leq n^{100}} \abs{g(\nu)} > \frac{w(n)}{\sqrt n} \abs{g(0)}?
\end{gather*}
\end{prob}
If we can solve Montgomery's problem for $B=100$ it is clear that there does not exist such a function  $w(x)$ in Tur\'an's problem. These conflicting guesses and a general lack of understanding of the situation has meant that until very recently  it has  not been clear to us what the correct order of magnitude in the problem should be. See the discussion in arXiv version 2 of \cite{Andersson3}.

\section{Main results}

The purpose of this paper is to solve Montgomery's problem and give the right order of magnitude for $(\star)$ if $m=\lfloor n^B \rfloor$. We will use an  estimate of character sums over finite fields of Katz. Before we state our theorems we state the following Lemma which we will prove in Section 3.
\begin{lem}
 Let $h \geq 2$ be an integer and let $q$ be a prime power. Then there exists unimodular complex numbers $z_1,\ldots z_q$ such that
 \begin{gather*}
   \max_{\nu=1,\ldots,q^h-2} \left|\sum_{k=1}^q z_k^\nu \right| \leq (h-1) \sqrt q. 
 \end{gather*}
\end{lem}

We first state a non explicit version of our result which has a somewhat easier proof than our sharpest version.
\begin{thm}
  Let $\delta>0$. Then we have uniformly for  $1+\delta \leq B \leq n$ the following estimate 
 $$   \sqrt{B n} \ll  \inf_{\abs{z_k} =1} \max_{\nu=1,\ldots,\lfloor n^B \rfloor} \abs{\sum_{k=1}^n  z_k^\nu} \ll B \sqrt n.  $$
  If $\abs{z_k} =1$ is replaced by $\abs{z_k} \geq 1$ then the result still holds for $2 \leq B \leq n$.
\end{thm}
This Theorem gives a direct solution to  Problem 1 of Montgomery and improves on Erd{\H o}s-Renyi's result Eq. \eqref{erd} for $2 < B \ll \log n$ by  removing a factor $\sqrt{(\log n)/B}$.
\begin{proof} The lower bound is Eq. \eqref{mont}. For the upper bound it is sufficient to consider $B=h$ for integers $h \geq 2$.
 Since $2^m$ is a prime power  we can choose $\{z_{k,m} \}_{k=1}^{2^m}$ in Lemma 1 such that 
\begin{gather*} \abs{\sum_{k=1}^{2^m} z_{k,m}^\nu} \leq \left(\lceil \frac N m \rceil-1 \right)2^{m/2}. \qquad (\nu=1,\ldots,2^N-2) \\ \intertext{If $a_m=0,1$ are the digits in the binary expansion of $n$, i.e.
$n=\sum_{m=0}^M a_m 2^m,$ then the sum $\sum_{m=0}^M a_m \sum_{k=1}^{2^m} z_{k,m}^\nu$ is a power sum of $n$ elements and can be written as $\sum_{k=1}^n z_k^\nu$. By the triangle inequality it is clear that}
 \abs{\sum_{k=1}^n z_k^\nu} \leq \sum_{m=0}^M   a_m \left(\lceil \frac N m \rceil-1 \right)2^{m/2} \ll \frac N M  2^{M/2} \\ \intertext{for $\nu=1,\ldots,2^N-2$. By the binary expansion of $n$ we have that   $n \leq 2^{M+1}-1$ and we see that  $n^h \leq 2^N-2$ when $n,h \geq 2$  and $N = h(M+1)$. We obtain the following estimate.}
 \abs{\sum_{k=1}^n z_k^\nu} \ll h \sqrt n.\qquad (\nu=1,\ldots,n^h)
\end{gather*}
\end{proof}

Theorem 1 implies the following Corollary.
\begin{cor}
 It is not possible to find a function $w(x)$ in Tur\'an's problem, Problem 2.
\end{cor}
\begin{proof} Choose $B=100$ in Theorem 1 and $b_1=\cdots=b_n=1$ in Tur\'an's problem. \end{proof}

We will now state a somewhat sharper version of Theorem 1, which we will
choose to state for $|z_k| \geq 1$ instead of $|z_k|=1$.
\begin{thm}
 Let $h \geq 2$ be an integer and $\epsilon>0$. One then has that
  \begin{gather*}
   C_h \sqrt n -O(n^{-1/2})
 \leq \inf_{\abs{z_k} \geq 1} \max_{\nu=1,\ldots,n^h} \abs{\sum_{k=1}^n
      z_k^\nu} \leq (h-1) \sqrt n+\Oh{n^{0.2625+\epsilon}}. 
  \end{gather*}
where $C_{2m}=(m!)^{1/2m} $ and $C_{2m+1}=C_{2m}$.
\end{thm}
We will prove this Theorem in Section 4.

\begin{rem} Theorem 2 still holds if $|z_k| \geq 1$ is replaced by  $\abs{z_k}= 1$ since the construction that yields the upper bound uses unimodular complex numbers. In this case the lower bound in Theorem 2 can be improved slightly by using the method from \cite{Andersson4}, e.g. for $h=3$ we get the lower bound $\sqrt {2n}$ instead of $\sqrt n$. At present there is however no $h > 2$ where we can attain asymptotic estimates in the problem. \end{rem}

By choosing $h=3$ in  Theorem  2 (with $|z_k|=1$) we get the following Corollary.
\begin{cor}
 Let $\alpha>1$ be a real number. One then has that
   \begin{gather*}
    \inf_{\abs{z_k}= 1} \max_{\nu=1,\ldots,\lfloor \alpha n^2 \rfloor } \abs{\sum_{k=1}^n
      z_k^\nu} \leq 2 \sqrt n+\Oh{n^{0.2625+\epsilon}}. \qquad (\epsilon>0)
  \end{gather*}
\end{cor}
This  improves on Theorems 6 and 7 from version 2 on arXiv of \cite{Andersson3}  for $\alpha > 4$. 
\begin{rem}
 In the proof of Theorem 2 we use the methods of \cite{Andersson3}, Section 6. The simpler methods of \cite{Andersson3}, Subsection 4.1 can be used to prove a weaker form of  Theorem 2 that still implies the Corollary.
\end{rem}
We will here suggest the following open problems (Compare with \cite{Andersson3}, Problem 3).
\begin{prob} 
  Find an increasing function $\Lambda(x)$ such that for each $B>1$ one has that
 \begin{gather*}
\inf_{\abs{z_k}= 1} \max_{\nu=1,\ldots,\lfloor n^B \rfloor } \abs{\sum_{k=1}^n
      z_k^\nu} \sim \Lambda(B) \sqrt n.
 \end{gather*}
\end{prob}
\begin{quest} Is the upper bound in Theorem 2 sharp? Can we choose $\Lambda(h)=h-1$ in Problem 3 when $h \geq 2$ is an integer? 
\end{quest}
Also the positive solution to Montgomery's problem suggests that the following problem might be possible to solve.
\begin{prob} 
 Let $h \geq 2$ be an integer. Find an increasing  function $\Lambda_h(\alpha)$ such that for each $\alpha>0$ 
 \begin{gather*}
\inf_{\abs{z_k}= 1} \max_{\nu=1,\ldots,\lfloor \alpha n^h \rfloor } \abs{\sum_{k=1}^n
      z_k^\nu} \sim \Lambda_h(\alpha) \sqrt n.
 \end{gather*}
\end{prob}
For $h=2$ this is given as Conjecture 4 in \cite{Andersson3}.
\begin{quest} Assume that we can solve Problem 4. Is it true that  $\lim_{\alpha \to \infty} \Lambda_h(\alpha) = \lim_{\alpha \to 0^+} \Lambda_{h+1}(\alpha)$?
\end{quest}
If Question 2 is true and Problem 3 solvable then we have that the function $\Lambda(x)$ in Problem 3 is piece wise constant for $n<x<n+1$.

\section{Proof of Lemma 1}

\begin{proof} Let $F$ be a finite field of order $q$ and let $E$ be an extension field of $F$ of order $q^h$. Let $x_1,\ldots,x_q$ denote the elements of $F$. Let $\omega$ be an element that generates the multiplicative group $E^*$, and let $\chi$ be a multiplicative character on $E$ of order $q^h-1$. Choose 
\begin{gather} \notag z_k=\chi(\omega+x_k). \qquad (k=1,\ldots,q) \intertext{Then} 
 \label{charsum}
 \sum_{k=1}^q z_k^\nu=\sum_{x \in F} \chi_\nu(\omega+x), \end{gather}
where $\chi_\nu=\chi^\nu$ is a non trivial character on $E$ unless $(q^h-1)|\nu$. By Theorem~1 of Katz \cite{Katz} its absolute value can be estimated from above by $(h-1)\sqrt n$ for $\nu=1,\ldots,q^h-2$.
\end{proof}
\begin{rem} The character sum Eq. \eqref{charsum} first occurred as eigenvalues of adjacency matrices of graphs in Chung \cite{Chung}. Katz proof uses Weil's estimates \cite{Weil}, although he prefers to use the language of Deligne.  \end{rem}

\begin{rem}
 For $h=2$ the construction in Lemma 1 coincides with our construction in Eq. (7) of \cite{Andersson2}, where we used a result of Bose \cite{Bose} on $B_2$ sequences $\pmod {q^2-1}$. Bose's result has been generalized to $B_h$  Sidon sequences  $\pmod{q^h-1}$ in Bose-Chowla \cite{BoseChowla}, and  this was how we originally arrived at the construction in Lemma 1. We will therefore describe  the construction of Bose-Chowla. Let $F$ and $E$ be a finite fields and $\omega$ an element in $E$ which is a multiplicative generator for $E^*$ as  as in the proof of Lemma 1. Choose the discrete logarithm so that $\log \omega=1$. Let  
\begin{gather}
  a_k=\log(\omega+x_k).  \qquad (k=1,\ldots,q)
\end{gather}
Then $\{a_k \}$ is a Sidon set of order $h$, or $B_h$ set $\pmod{q^h-1}$, or in other words
the equation $a_{j_1}+\cdots+a_{j_h} \equiv a_{i_1}+\cdots+a_{i_h} \pmod {q^h-1}$ has no non-trivial solution.  Since the case $h=2$ had been  successfully used in the problem it was natural to expect that this construction could give good estimates for a general $h \geq 3$ and  $\nu=1,\ldots,p^h-2$. Indeed, numerical investigation (with $h=3,4$ and small primes)     suggested to us that if 
\begin{gather*}
  z_k=e\left( \frac{a_k}{q^h-1} \right), \qquad (k=1,\ldots,q)
\end{gather*}
then  we should have 
\begin{gather*}
 \abs{\sum_{k=1}^q z_k^\nu} \leq (h-1) \sqrt n. \qquad (\nu=1,\ldots,q^h-2)
\end{gather*}
This can explained as follows: We can define the character $\chi$  by
\begin{gather*}
 \chi(\omega^n)= e \left( \frac n {q^h-1} \right), \qquad \chi(0)=0. \\ \intertext{Then}
 z_k=\chi(\omega+x_k),
\end{gather*}
where $\chi$ is a character of order $q^h-1$ on  $E^*$ as in the proof of Lemma 1.
\end{rem}

\section{Proof of Theorem 2}

The lower estimate follows from Theorem 2 in our paper \cite{Andersson1}
\begin{gather*}
 \left( \frac {n! m!}{(n-m)!} \right)^{1/(2m)} \leq \max_{\nu=1,\ldots,n^{2m}} \abs{\sum_{k=1}^n z_k^\nu}. \qquad (|z_k|\geq 1)
\end{gather*}

Hence we will concentrate in proving the upper bound. We will use the method from \cite{Andersson3}. We have the following Lemma:

\begin{lem} \label{fund}   Let $\epsilon>0$, $0<\theta<1$,   $C \geq 1$ and let $h \geq 2$ be an integer. Suppose that  $(z_1,\ldots,z_n)$ is an $n-$tuple of  unimodular complex numbers,
 \begin{gather} \label{star2}
  m \sim n^\theta, \\ \intertext{and} \notag
    \abs{\sum_{k=1}^n z_k^\nu} \leq C    \sqrt n. \qquad \qquad (\nu=1,\ldots, n^h) \\ \intertext{Let $\mathcal N=\{1,\ldots,n \}$. Then there exist a subset $\mathcal M_0  \subset  \mathcal N$, with $\# \mathcal M_0=m$  such that }
    \abs{\sum_{k \in \mathcal M_0} z_{k}^\nu} \ll_{\epsilon} m^{1/2+\epsilon}. \qquad  \qquad(\nu=1,\ldots, n^h) \notag
  \end{gather}
\end{lem}
\begin{proof} The Lemma is the same as Lemma 9 from \cite{Andersson3} except for the fact that  $\lfloor \alpha n^2 \rfloor$ is replaced by $n^h$. 
The proof is the same, but we need the choice $N >\frac h {2 \theta \epsilon}$ instead of $N>\frac 1 {\theta \epsilon}$. \end{proof}

 By the Baker-Harman-Pintz theorem \cite{BakerHarmanPintz} we can choose a prime $n<p$ such that $p-n \asymp p^{0.525}$. 
By the construction given in Lemma 1, we can choose a $p-$tuple $(z_1,\ldots,z_{p})$   of unimodular complex numbers such that
\begin{gather*}
  \abs{\sum_{k=1}^{p} z_k^\nu} \leq (h-1) \sqrt {p}. \qquad  \qquad (\nu=1,\ldots,p^h-2) 
\end{gather*}
Let $m=p-n$. By  Lemma \ref{fund} with $\theta=0.525$ and $C=h-1$ we can choose a subset $\mathcal M_0 \subset \{1,\ldots,p\}$ with $\# M_0=m$ such that
\begin{gather*}
  \abs{\sum_{k \in \mathcal M_0} z_k^\nu} \leq  n^{0.2625+\epsilon}. \qquad \qquad (\nu=1,\ldots,p^h-2) \qquad (\epsilon>0)
\end{gather*}
Let $\mathcal N = \{1,\ldots,p \} \setminus \mathcal M_0$. It is clear that $\# \mathcal N=n$ and by the triangle inequality it follows for $1 \leq \nu \leq n^h \leq p^h-2$ that
\begin{gather*}
     \begin{split}
     \abs{\sum_{k \in \mathcal N} z_k^\nu} &=  \abs{\sum_{k=1}^{p} z_k^\nu - \sum_{k \in \mathcal M_0} z_k^\nu}, \\ 
        &=  \abs{\sum_{k=1}^{p} z_k^\nu}+ \Oh{\abs{\sum_{k \in \mathcal M_0} z_k^\nu}},  \\
        &\leq (h-1) \sqrt{n+\Oh{n^{0.525}}} + \Oh{n^{0.2625+\epsilon}}, \\ &\leq
            (h-1) \sqrt n + \Oh{n^{0.2625+\epsilon}}.
   \end{split}
\end{gather*}
which finishes the proof of our Theorem. \qed

\bibliographystyle{plain}

\end{document}